\documentclass[british,smallextended]{svjour3}
\usepackage[T1]{fontenc}
\usepackage[latin9]{inputenc}
\usepackage{babel}
\usepackage{array}
\usepackage{float}
\usepackage{textcomp}
\usepackage{url}
\usepackage{amstext}
\usepackage{graphicx}
\usepackage[unicode=true]
 {hyperref}

\makeatletter

\providecommand{\tabularnewline}{\\}

\RequirePackage{fix-cm}

\smartqed  

\@ifundefined{showcaptionsetup}{}{%
 \PassOptionsToPackage{caption=false}{subfig}}
\usepackage{subfig}
\makeatother

\begin{document}

\title{Comparison Between Robust and Stochastic Optimisation for Long-term
Reservoir Management Under Uncertainty}

\author{Thibaut Cuvelier\and Pierre Archambeau\and Benjamin Dewals\and Quentin
Louveaux}

\institute{Thibaut Cuvelier \at \\
\email{tcuvelier@ulg.ac.be}}

\date{Received: date / Accepted: date}
\maketitle
\begin{abstract}
Long-term reservoir management often uses bounds on the reservoir
level, between which the operator can work. However, these bounds
are not always kept up-to-date with the latest knowledge about the
reservoir drainage area, and thus become obsolete. The main difficulty
with bounds computation is to correctly take into account the high
uncertainty about the inflow to the reservoir. We propose a methodology
to derive minimum bounds while providing formal guarantees about the
quality of the obtained solutions. The uncertainty is embedded using
either stochastic or robust programming; we derived linear models,
which ensures good performance to compute the solutions.  
\end{abstract}

\section{Introduction\label{sec:Introduction}}

Drinking-water production is a vital need, and tap water is a basic
service that may not be interrupted. This water may come from multiple
sources, such as groundwater or dams (surface water). Utility managers
cannot afford lacking water to inject in their distribution system,
and this is one reason that impelled them to build large dams and
reservoirs. Those serve also multiple other purposes, such as hydropower
\cite{Bieri2013} or flood control \cite{Camnasio2011}. 

These reservoirs must deliver the expected level of service with a
very low probability of failure. A common management technique is
to use predefined rules, i.e. \emph{rule curves}, indicating upper
and lower bounds on the level of the reservoir at any period of time
throughout the year \cite{Chang2005}. Modern operation systems usually
prefer real-time control \cite{Schwanenberg2015}, with no predefined
rules, but this is harder to implement than traditional rule curves. 

Rule-curve-based management is expected to take full benefit of the
potential of the reservoirs. However, the rule curves must be periodically
updated; otherwise, they become less reliable over time and may eventually
fail to provide the expected service. The quality degradation of these
operational rules is mostly due to external evolution\textemdash be
they in the required amount of water or in climate change. 

In this context, \emph{formal} optimality guarantees of the rule curve
are an asset for the operator: if there is any lack of water, the
operator is able to prove that every possible action was taken to
prevent this situation. To this end, mathematical optimisation can
be used to derive the rule curves; with a well-defined and reproducible
computational framework, the rule curves can also be regularly updated,
and thus remain relevant, even under external evolutions. Most current
optimisation techniques \cite{Maier2014} are based on evolutionary
computations \cite{Ahmad2014} (mainly genetic algorithms \cite{Chang2005},
sometimes coupled with simulation \cite{Taghian2014}; newcomers like
harmony search are gaining traction in the community \cite{Bashiri-Atrabi2015}).
In contrast to them, mathematical optimisation \cite{Labadie2004,Zhang2015}
can provide global-optimality guarantees while having acceptable computational
time requirements (under some assumptions, such as convexity), which
this article demonstrates in the case of long-term rule curves. 

Whatever the optimisation method, it must take into account the inherent
uncertainty in the input data (namely, the inflow). In mathematical
optimisation, there are two well-known paradigms to handle uncertainty:
stochastic and robust optimisation. 
\begin{enumerate}
\item The first is the most common one \cite{WRCR:WRCR3909}, and is also
known as \emph{explicit stochastic optimisation} \cite{Celeste2009}.
It implies that the model directly uses probabilistic information
(usually, by sampling the probability density function, i.e. using
inflow scenarios for each river). The objective function minimises
some risk measure \cite{Shapiro2014}, often the expected value of
the deterministic function. However, in our case, the maximum function
makes more sense than the expectation (see Section \ref{subsec:Robust-stochastic-model}
for a full discussion). 
\item Recently, another approach has been explored: \emph{robust optimisation}
\cite{Pan2015}, whose roots lie in mathematical optimisation communities
\cite{Ben-Tal2009}. It considers the uncertain values as pertaining
to a so-called \emph{uncertainty set} and it optimises for the worst
case within that set; a common choice is to use the confidence intervals
around the average value. It can be considered as a kind of \emph{implicit}
stochastic optimisation (ISO) \cite{Celeste2009}, as the actual model
is deterministic, albeit working with data that is adapted to the
uncertainty (without any need for hedging the solutions to different
sets of possible inflow scenarios, as is usually necessary in ISO
approaches, like in \cite{Zhao2014}). 
\end{enumerate}
In this paper, we compare these two paradigms based on a real-world
case study of a Belgian dam on the river Vesdre, at Eupen. Historically,
a very conservative process was used to compute the rule curves, with
a single objective: water supply. We discuss to which extent the two
considered approaches, robust and stochastic optimisation, may contribute
to give more freedom to the operator for other purposes, such as hydropower. 

The article is structured as follows. After a brief description of
the case study in Section \ref{sec:Case-study}, the optimisation
models corresponding to each uncertainty paradigm are detailed in
Section \ref{sec:Methodology}. The results of the optimisation process
are discussed in Section \ref{sec:Results}. Finally, conclusions
are drawn in Section \ref{sec:Conclusion}, with further directions
to improve this work. 

\section{Case study \label{sec:Case-study}}

In this study, we consider the Eupen dam (also called the Vesdre dam),
which lies on the river Vesdre in Eastern Belgium. The reservoir is
fed by its natural drainage area (6920 ha) and by another river (the
Helle), from which a diversion tunnel was built. This diversion increases
the effective catchment area to 10,595 ha. In normal-operation mode,
the tunnel is open and only a minimum environmental flow remains in
the river Helle. If the reservoir level reaches its maximum level,
the tunnel can be closed, so that the whole discharge of the river
Helle is conveyed through its natural riverbed. The dam height is
50 m and the maximum storage capacity is 25 hm\textthreesuperior .
Currently, the main purpose of the Vesdre reservoir is to provide
drinking water to its surroundings. It can also be used to produce
hydropower (2.6 MW) and to help controlling floods and low flows. 

For now, its operation uses empirical rules, which are functions of
the measured discharges and the weather forecasts, e.g. in order to
prevent flooding at some points in the hydraulic network. Drinking-water
availability is guaranteed through a minimum \emph{rule curve}: the
reservoir level may not drop below a certain threshold (which varies
over the year). 

This curve was computed once, in the 1980s, using a very conservative
process. A worst-case scenario was derived from historical data: each
month in this scenario was taken as the driest one in history; then,
the lowest acceptable level was determined to ensure the expected-water
supply guarantee. In the following, we analyse how stochastic and
robust optimisation may contribute to revise this rule curve, by making
use of today's enhanced computational power and more recent data. 

We use as main input data the characteristics of the catchment, reservoir,
and dam (see Table \ref{tab:Technical-characteristics} in appendix),
as well as time series of recorded natural inflow discharge to the
reservoir and flow rates in the Helle river during 24 calendar years
(1992-2015). 

\section{Methodology \label{sec:Methodology}}

We develop optimisation models to compute the minimum acceptable level
in a reservoir to ensure that drinking water availability will not
be exposed. More specifically, starting from the beginning of a drought,
the available drinking water storage must be sufficient to enable
water supply during a predefined period of time. In the considered
case study, the risk level is set such that water production for two
years is ensured, but the methodology does not depend on the specific
value of this parameter. 

\subsection{Uncertain reservoir models \label{subsec:Uncertain-models}}
\begin{remark}
In the following mathematical notations, the optimisation variables
are indicated by a \textbf{bold} font. In other words, the value for
these symbols is the result of the optimisation process, and are not
fixed beforehand. 
\end{remark}

\subsubsection{Deterministic model \label{subsec:Deterministic-model}}

The real-world decision variable is the \emph{reservoir level}, but
it is not directly used in the model: we use the \emph{volume} of
water instead, which is in one-to-one mapping to the reservoir level.
This variable is denoted by $\mathbf{storage}_{t}$ at time $t$ ($\mathrm{m}^{3}$).
The optimisation problem consists mainly of one equation, the mass
balance over a time step \cite{WRCR:WRCR3909,Arunkumar2012,Pan2015}:
\[
\mathbf{storage}_{t+1}=\mathbf{storage}_{t}-\mathbf{output}_{t}+\mathbf{input}_{t},\qquad\forall t,
\]
where the inputs correspond to the inflowing rivers, and the outputs
to the various dam purposes (drinking water, hydropower, etc.). \\
The inputs are composed of tributary rivers (the set $\mathit{tributaries}$)
and of diverted rivers (the set $\mathit{diverted}$): 
\[
\mathbf{input}_{t}=\sum_{r\in\mathit{tributaries}}\mathrm{flow}_{t,r}+\sum_{r\in\mathit{diverted}}\mathbf{diverted}_{t,r},\qquad\forall t.
\]
For the considered case, the output from the reservoir is the sum
of the drinking water consumption, the minimum environmental flow
in the river, and the release that is not useful to fulfil the reservoir's
purposes: 
\[
\mathbf{output}_{t}=\underbrace{\mathrm{drinkingWater}_{t}+\mathrm{environmentalFlow}_{t}}_{\mathrm{constant\,(i.e.\,not\,fixed\,by\,the\,optimisation\,process)}}+\,\mathbf{release}_{t},\qquad\forall t.
\]
In this equation, several terms are fixed by the operator before the
optimisation (drinking water and environmental flow requirements);
the only actual decision variables are the releases. Other losses
could be taken into account, such as evaporation, but they are negligible
in the considered area \cite{Finch2008}. 

The reservoir level may not drop below a minimum threshold, defined
such that the operator can still extract water from the reservoir
for water supply (it corresponds to the position of water outlets,
and not to any risk level). Similarly, it has a maximum level that
is fixed beforehand to keep a safety margin in case of floods. These
constraints can be expressed as: 
\[
\mathrm{minStorage}\leq\mathbf{storage}_{t}\leq\mathrm{maxStorage},\qquad\forall t.
\]
Likewise, the contributions of the diverted rivers can be decided
up to some level, as long as two constraints are respected: a minimum
environmental flow must remain in the diverted river (and thus cannot
flow into the reservoir), while the maximum discharge capacity through
the diversion pipe (due to its dimensions) constitutes an upper bound.
\[
\mathbf{diverted}_{t,r}\leq\mathrm{maxDischarge}_{r},\qquad\forall t,\forall r\in\mathrm{diverted}.
\]
\[
\mathbf{diverted}_{t,r}\leq\mathrm{flow}_{t,r}-\mathrm{environmentalFlow}_{r},\qquad\forall t,\forall r\in\mathrm{diverted}.
\]
Also, the release from the dam $\mathbf{release}_{t}$ has an upper
bound, related to the hydraulic capacity of the hydropower plant and
the bottom outlets. 
\[
\mathbf{release}_{t}\leq\mathrm{penstockHydropower}+\mathrm{bottomOutlet,\qquad}\forall t.
\]

Due to the relative position between the spillway crest and the maximum
allowable reservoir level, the spillway is not taken into account
in this model. 

The goal is to determine an enhanced rule curve, i.e. a new lower
bound for the reservoir level: the objective function minimises the
total stored volume throughout the year, with all periods having the
same weight: 
\[
\min\sum_{t}\mathbf{storage}_{t}.
\]
This way, the solution is the most critical situation while still
being feasible from the beginning to the end of the optimisation horizon. 
\begin{remark}
This objective function is supported by the fact that, if the solver
lowers the value for one time step at the expense of another, the
variations for these two time steps have the same effect on the objective
value. In other words, at the optimality, it is not possible to decrease
the total amount of water that is stored throughout the year, and
changing the obtained value for any time step forces to reconsider
the solution at other time steps. 
\end{remark}

All constraints detailed above and the objective are linear: this
optimisation problem is thus very tractable (large instances can be
solved quickly) \cite{Vanderbei2014}. However, the model ignores
the uncertainty in the inflow (i.e. $\mathrm{flow}_{t,r}$): it can
only consider one inflow scenario, and optimises over that scenario,
which is not necessarily representative of basin's dynamics. The following
sections present two approaches to incorporate the uncertainty into
the model while keeping it linear. 

\subsubsection{Stochastic model \label{subsec:Robust-stochastic-model}}

A first approach considers the inflow as stochastic, the uncertainty
being modelled as a series of scenarios \cite{Shapiro2014}; the rule
curve is such that the drinking-water requirement is guaranteed, whatever
the inflow scenario is. To optimise the rule curve, we simulate each
scenario independently; then, the actual solution to the optimisation
problem is the \emph{upper envelope} of the solutions to the individual
scenarios (a convex \emph{risk measure}, in \cite{Shapiro2014}).
This risk-averse modelling ensures that the computed rule curve is
conservative enough, based on the information that is known. 

Let $\mathbf{storage}_{t}^{s}$ denote the solution at time $t$ for
scenario $s$, and $\mathbf{ruleStorage}_{t}$ the actual value for
the rule curve at time $t$. The rule curve must be \emph{above} the
minimum level for any scenario (as it is the maximum of all the solutions),
which is translated by the constraint $\mathbf{storage}_{t}^{s}\leq\mathbf{ruleStorage}_{t}$.
The complete model is thus:
\[
\begin{array}{clc}
\min & \sum_{t}\mathbf{ruleStorage}_{t}\\
\mbox{s.t.} & \mathbf{storage}_{t}^{s}\leq\mathbf{ruleStorage}_{t} & \forall t,\forall s,\\
 & \mathbf{storage}_{t+1}^{s}=\mathbf{storage}_{t}^{s}-\mathbf{output}_{t}^{s}+\mathbf{input}_{t}^{s} & \forall t,\forall s,\\
 & \mathbf{input}_{t}^{s}=\sum_{r\in\mathit{tributaries}}\mathrm{flow}_{t,r}^{s}+\sum_{r\in\mathit{diverted}}\mathbf{diverted}_{t,r}^{s} & \forall t,\forall s,\\
 & \mathbf{output}_{t}^{s}=\mathrm{drinkingWater}_{t}+\mathrm{environmentalFlow}_{t}+\mathbf{release}_{t}^{s} & \forall t,\forall s,\\
 & \mathrm{minStorage}\leq\mathbf{storage}_{t}^{s}\leq\mathrm{maxStorage} & \forall t,\forall s,\\
 & \mathbf{diverted}_{t,r}^{s}\leq\mathrm{maxDischarge}_{r} & \forall t,\forall s,\forall r\in\mathit{diverted},\\
 & \mathbf{diverted}_{t,r}^{s}\leq\mathrm{flow}_{t,r}^{s}-\mathrm{environmentalFlow}_{r} & \forall t,\forall s,\forall r\in\mathit{diverted},\\
 & \mathbf{release}_{t}^{s}\leq\mathrm{penstockHydropower}+\mathrm{bottomOutlet} & \forall t,\forall s,\\
 & \mathbf{ruleStorage}_{t}\geq0 & \forall t,\\
 & \mathbf{storage}_{t}^{s}\geq0,\qquad\mathbf{output}_{t}^{s}\geq0,\\
 & \mathbf{input}_{t}^{s}\geq0,\qquad\mathbf{release}_{t}^{s}\geq0 & \forall t,\forall s,\\
 & \mathbf{diverted}_{t,r}^{s}\geq0 & \forall t,\forall s,\forall r\in\mathit{diverted},
\end{array}
\]

\paragraph{\textbf{Scenario generation. }}

An interesting point is the way of generating the scenarios. As a
property must be guaranteed for two years, each scenario should be
two-year long. They can correspond to two successive years in the
historical data (``merging''); another way of generating them is
to consider all possible pairs of one-year scenarios (``mixing'').
The latter technique destroys any kind of inter-year correlation (as
opposed to the first one), even though both keep intra-year correlations. 

Traditional stochastic models are multistage \cite{Labadie2004},
i.e. they consider a point in time where the operator is allowed to
reconsider their decisions (this includes stochastic dynamic programming).
However, this formalism would be detrimental to our use case: our
goal is to determine a rule curve that does not depend on the actual
scenario, thereby giving the operator a sure lower bound for their
reservoir management. This case is more pessimistic than what multistage
stochastic optimisation allows: our solution must be valid in all
cases at all time steps, while a multistage formulation considers
one solution per branch in the scenario tree \cite{Birge2011}; as
a consequence, the present approach is more suited for strategic planning
of the management, while multistage models apply for real-time decisions. 

\subsubsection{Robust model \label{subsec:Robust-model}}

The second uncertain model considers the inflow as belonging to an
uncertainty set, which we chose to be the confidence interval of the
inflow at the corresponding time period based on the historical data.
This choice is sometimes called \emph{interval uncertainty} or \emph{Soyster's
model} \cite{Ben-Tal2002}. 

The worst inflow to the reservoir can be computed explicitly: it corresponds
to the minimum inflow within the confidence interval (i.e. the lower
bound). As such, the robust model is very similar to the basic deterministic
one that is detailed in Section \ref{subsec:Deterministic-model},
except in the way the inflow is chosen: it is no more raw historical
data, but rather this worst case, computed with the historical data. 

\paragraph{\textbf{Confidence intervals. }}

Multiple ways of computing those confidence intervals (and thus their
lower bounds) can be thought of. The simplest approach is probably
to separately consider periods of time, and applying a standard statistical
technique (e.g., based on t-Student distribution), albeit this has
limited physical meaning. Other methods include fitting a dedicated
statistical model (such as the one in \cite{Adam2014}) and using
standard time series analysis procedures (like ARMA \cite{Abrahart1998}). 

\subsection{Model predictive control (MPC) \label{subsec:MPC}}

Uncertain models as such have one important defect: they guarantee
the two-year supply only for the first time step of the solution;
for the following ones, the guarantee is limited by the time horizon
of the model (i.e. the second time step has a guarantee for two years
\emph{minus} one time step). Using longer scenarios would not solve
this issue: for example, with three-year scenarios, the solution for
the first time step would be guaranteed for three years, which is
more constraining than expected. 

\begin{figure}
\begin{centering}
\includegraphics[scale=0.7]{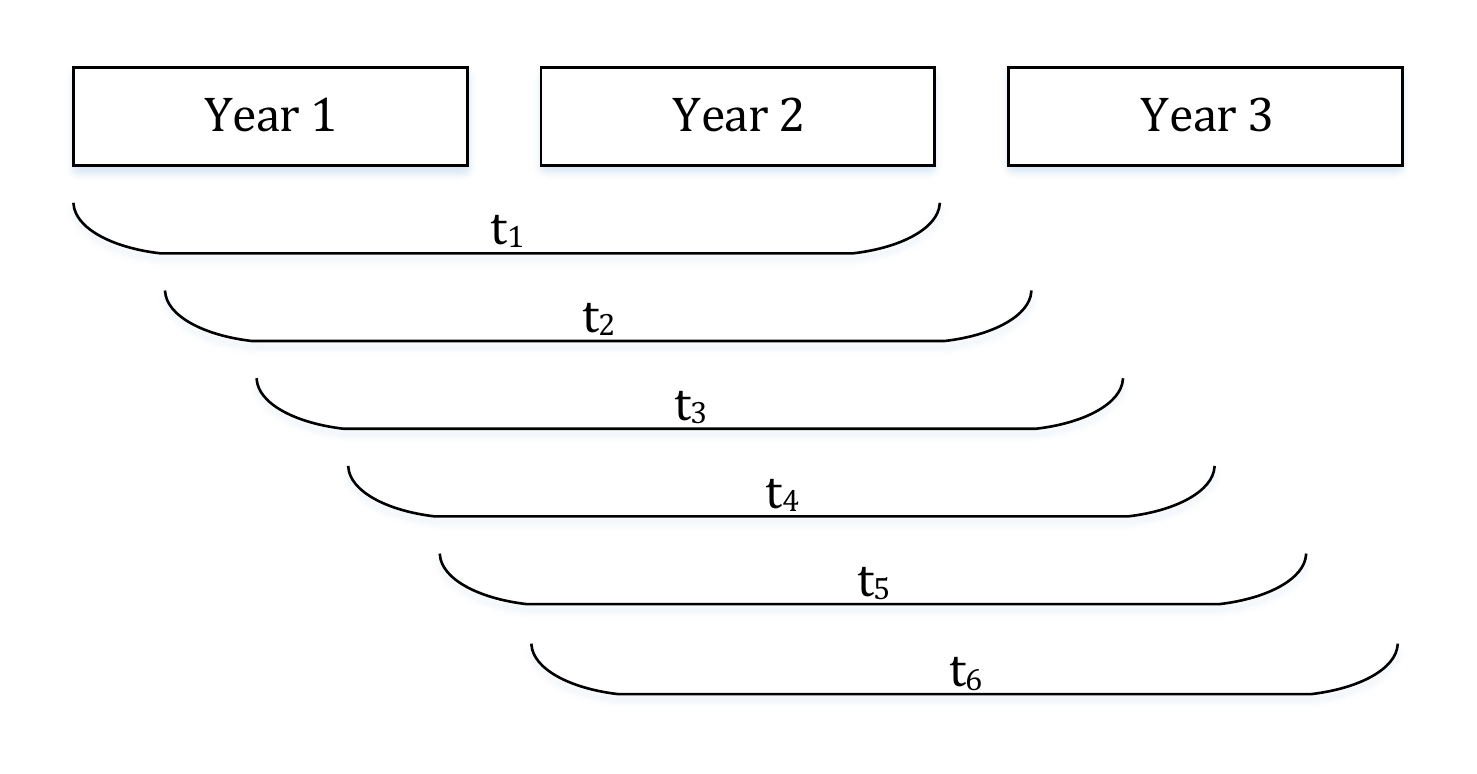}
\par\end{centering}
\caption{\label{fig:Behaviour-of-MPC}Behaviour of the MPC algorithm. The solution
for the first time step $t_{1}$ has an optimisation horizon limited
to the first two years (i.e. ensures the drinking-water guarantee
for two years); the one for $t_{2}$ is computed for a two-year period
starting at the second time step, i.e. the optimisation horizon of
$t_{1}$ shifted by one time step. }
\end{figure}

To work around this deficiency, we use moving time horizons: the algorithm
uses three-year scenarios, but the actual optimisation (performed
with an uncertain model, as explained earlier) merely happens on two
consecutive years of the scenario; only the solution for the first
time step is used to build the final solution (with the two-year guarantees).
To construct the complete solution (for one year), this two-year horizon
shifts, as depicted in Figure~\ref{fig:Behaviour-of-MPC}. 

This technique is called \emph{receding horizon control} \cite{Kwon2005},
and is often used in industrial process control, including beforehand
operational-rules computation (which is precisely the case here).
It is part of a more generic framework, \emph{model predictive control}
(MPC), that has already been used with great success in water resources
management, but mostly for real-time control \cite{Talsma2013,Becker2014}. 

As a consequence, to fully ensure the water supply guarantees, the
uncertain models must be solved a high number of times\textemdash albeit
these are completely independent computations and can be performed
in parallel. 

Thanks to this algorithm, all obtained solutions are year-to-year
continuous: there is no gap in the computed rule curve between the
end of a year and the beginning of the next one. The solutions of
the basic uncertain models of Section \ref{subsec:Uncertain-models}
do not have this continuity property, which is important for operators
to implement the rule curve. 

\section{Results \label{sec:Results}}

\subsection{Outcome of the optimisation}

\begin{figure}
\begin{centering}
\includegraphics[scale=0.4]{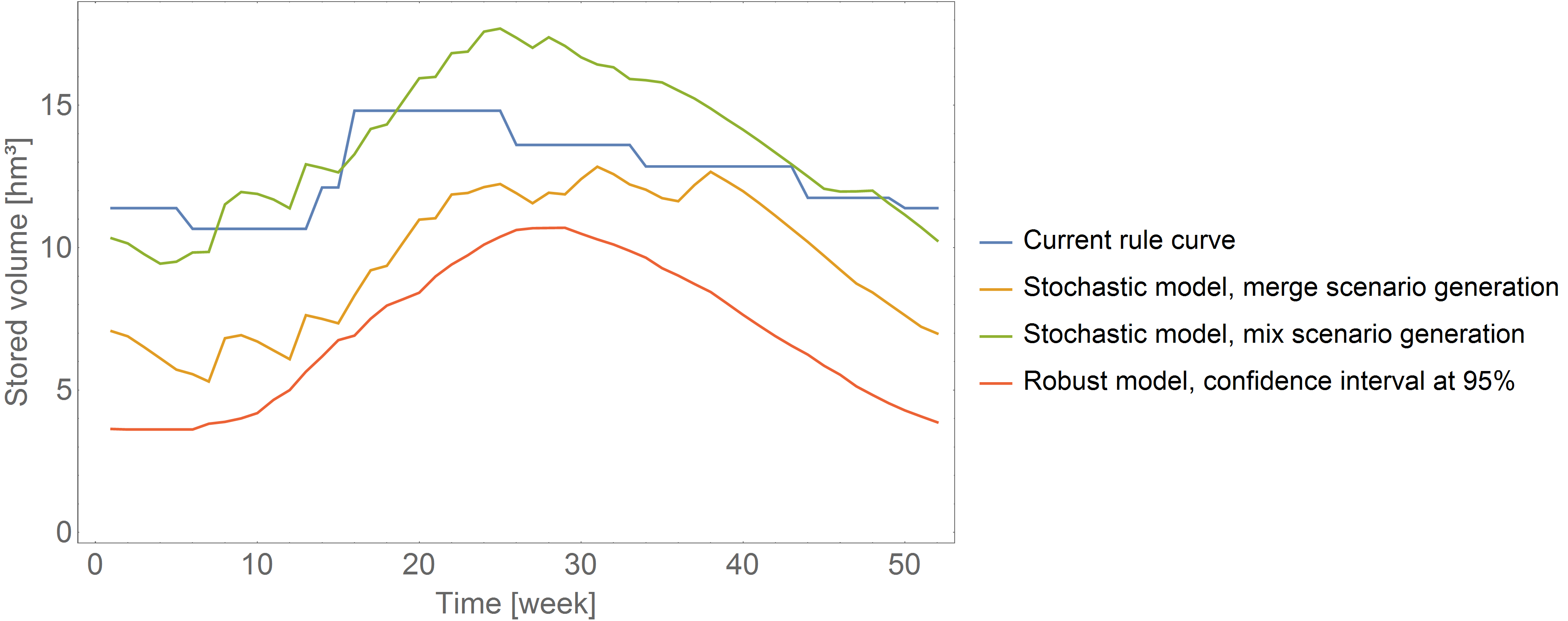}
\par\end{centering}
\caption{\label{fig:Compare-rule-curves}Comparing the two uncertainty models
to the current rule curve: on the one hand, stochastic, with two scenario
generation techniques; on the other one, robust.}
\end{figure}

The two described models (stochastic and robust models, both with
model predictive control) have been implemented in a Julia package
\cite{Bezanson2014}, \emph{\href{https://github.com/dourouc05/ReservoirManagement.jl}{ReservoirManagement.jl (https://github.com/dourouc05/ReservoirManagement.jl)}},
which is freely available on GitHub; it is based on the JuMP mathematical
modelling layer \emph{\cite{Dunning2015}}. Both models are compared
to the existing minimum rule curve in Figure~\ref{fig:Compare-rule-curves}.
Multiple such curves have been computed for each uncertainty model: 
\begin{itemize}
\item for the stochastic model: the two scenario generation techniques (merging
and mixing); 
\item for the robust model: multiple confidence intervals (based on a $t$-Student
model), starting at $95\%$ and up to $98.5\%$, with increasing conservativeness
(as shown in appendix, Figure~\ref{fig:Compare-rule-curves-confidence-level}).
Higher confidence levels ($99\%$ and beyond) do not allow for any
solution: too little water is available in the corresponding scenario. 
\end{itemize}
Overall, the computed rule curves strongly depend on how the uncertainty
is handled, i.e. stochastic or robust approach; nonetheless, the curves
remain relatively close to the current rule curve; yet, at some time
steps, all the proposed models are below the current rule curve. In
other words, depending on the way to model the uncertainty, the current
rule curve is either too conservative (none of our solutions needs
a reservoir level as high as prescribed by the existing rule curve)
or marginally unsafe (one model proposes to keep a slightly higher
level for about one third of the year). 

\subsection{Associating confidence levels}

\begin{figure}
\begin{centering}
\includegraphics[scale=0.4]{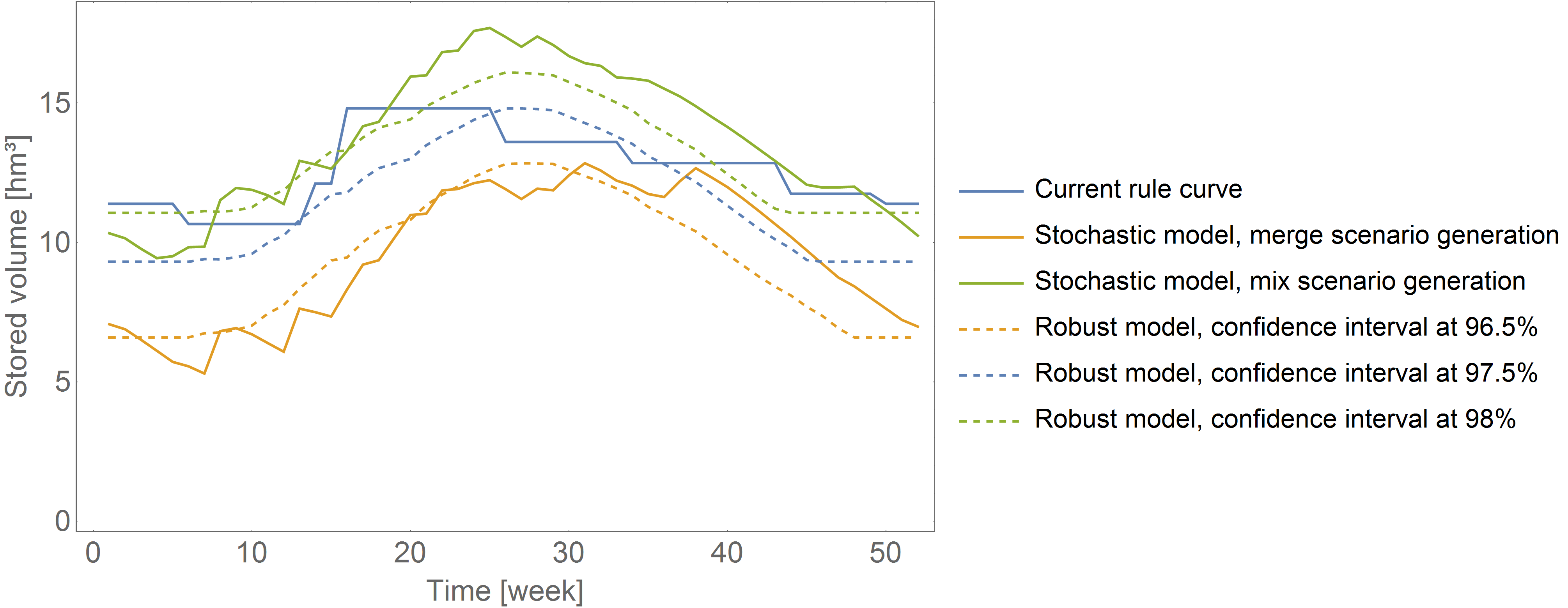}
\par\end{centering}
\caption{\label{fig:Compare-stochastic-robust-confidence-levels}Closest confidence
intervals corresponding to the stochastic solutions.}
\end{figure}

Based on both the stochastic and the robust approaches, we can estimate
a confidence level of feasibility for the stochastic solutions, by
choosing the robust solution that is the closest to the stochastic
curve to analyse within a fixed set of confidence levels (between
$95\%$ and $99\%$). This correspondence is shown in Figure~\ref{fig:Compare-stochastic-robust-confidence-levels}.
We can say that the merging-scenario-generation technique gives a
confidence level of approximately $96.5\%$, the existing rule curve
$97.5\%$, the mixing scenario generation $98\%$. 

\subsection{Most important scenarios for stochastic optimisation}

\begin{figure}
\subfloat[\label{fig:Impact-scenario-on-solution-yearly}The colour of the scenario
curve corresponds to the \emph{average discharge}: the bluer curves
indicate a wetter year.]{\begin{centering}
\includegraphics[scale=0.4]{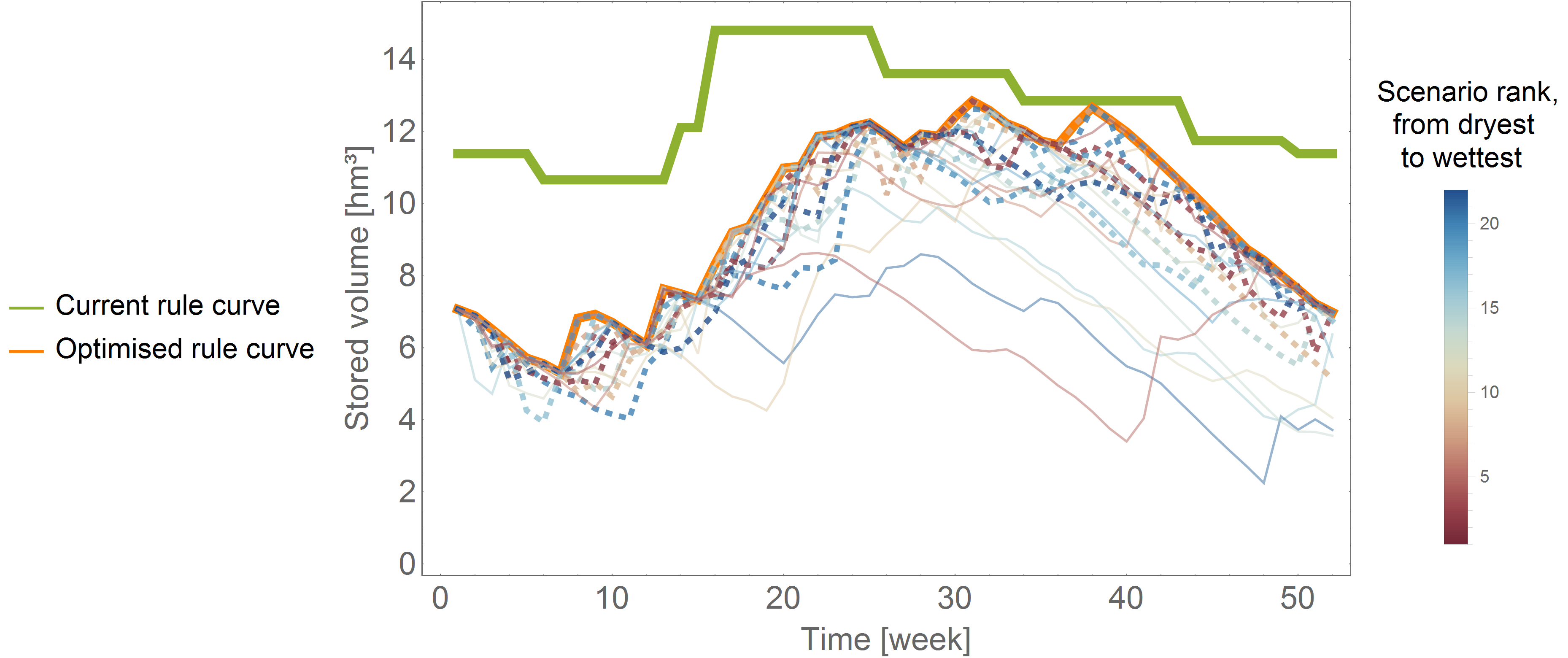}
\par\end{centering}
}
\begin{centering}
\subfloat[\label{fig:Impact-scenario-on-solution-driest-month}The colours correspond
to the average discharge during \emph{the driest month}: the reddest
curve indicates the driest month among all the scenarios. ]{\begin{centering}
\includegraphics[scale=0.4]{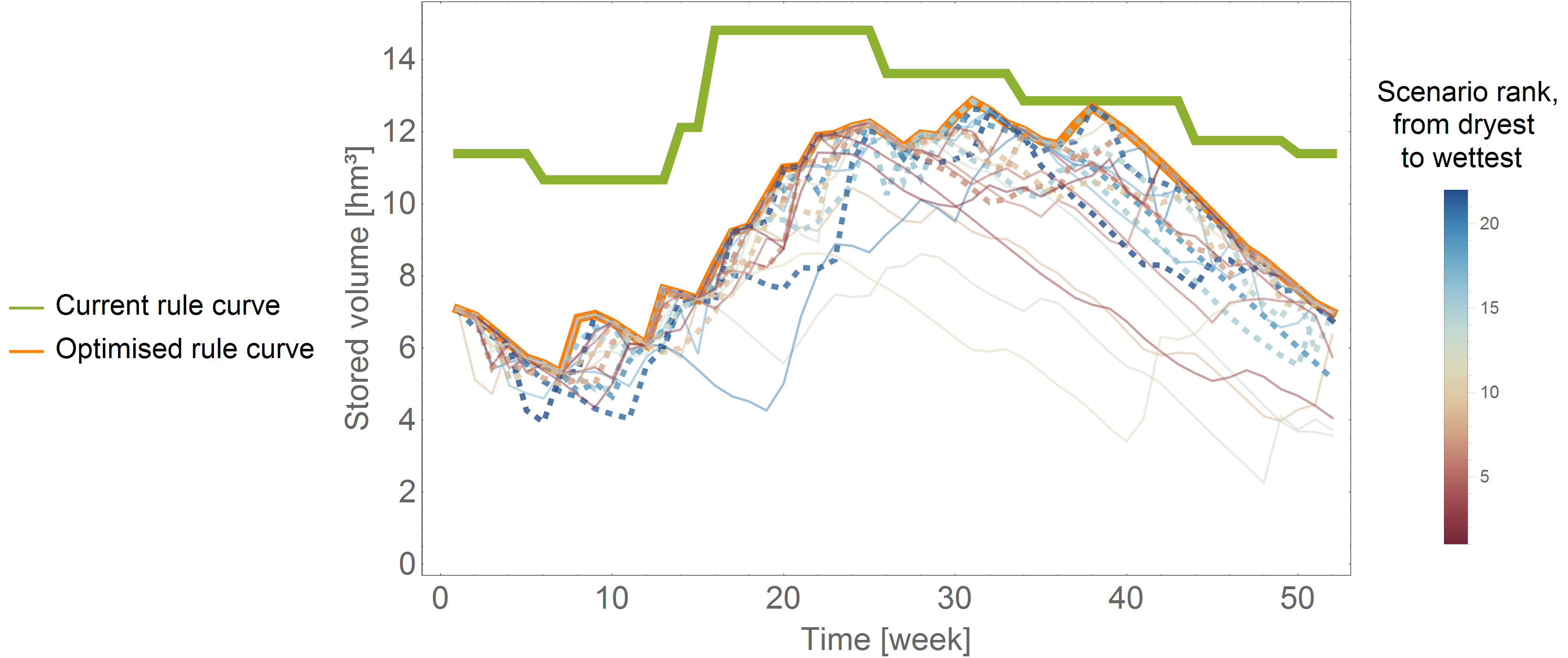}
\par\end{centering}
}
\par\end{centering}
\caption{Impact of the scenarios on the solution for a stochastic solver (here,
merging is depicted, without MPC). The optimisation curves (i.e. the
result of the optimisation and the scenarios) correspond to a 2-year
scenario; each of them is plotted against the original rule curve
for the first year. The colour scales are relative: dark blue indicates
the wettest scenario, dark-red the driest scenario. A dashed line
indicates a support scenario. }
\end{figure}

For the stochastic model (Section \ref{subsec:Robust-stochastic-model}),
there is not a single scenario that fully defines a given solution,
as shown in Figure~\ref{fig:Impact-scenario-on-solution-yearly}:
instead, a limited number of scenarios have an impact on the rule
curve (these scenarios may be called \emph{support vectors} \cite{Cortes1995});
the other scenarios have no influence on the solution. 

Among them, some are closer to the wettest year, others to the driest
one (Figure~\ref{fig:Impact-scenario-on-solution-yearly}): a low
average discharge throughout the scenario does not imply that it yields
the most conservative solution; the distribution of the inflow over
the year has a substantial influence on the result. This means that
limiting the study to the driest years would not be enough to derive
a sufficiently reliable minimum rule curve. In contrast, a more important
factor is the driest month (as depicted in Figure~\ref{fig:Impact-scenario-on-solution-driest-month}),
as the support scenarios correspond instead to those containing some
of the driest months, and these define the solution for a large period
of time. Figures \ref{fig:Impact-scenario-on-solution-driest-three-months}
to \ref{fig:Impact-scenario-on-solution-wet} in appendix present
similar results for averaging periods of three or six months, and
of the wet and dry seasons. The highest correlation is seen for the
driest month, but the three driest months and the dry season are similarly
correlated, as indicated in Table~\ref{tab:Impact-scenario-on-solution}.
The support scenarios are not only made up of dry years, but also
\emph{a few} very wet ones. Also, these results suggest that the month
is a very relevant time scale for defining the rule curve. 

The main conclusion for the optimisation of this analysis is that
the time step for the optimisation of this model should never be longer
than one month, otherwise the resulting rule curve may miss some important
events. It also gives an easy (but approximate) criterion to discriminate
support scenarios from redundant ones, based on the driest month contained
in the time series. 

\begin{table}
\caption{\label{tab:Impact-scenario-on-solution} Analysis of the time periods
that define the optimisation results. There are 22 scenarios, of which
9 are support scenarios. }

\begin{tabular}{>{\raggedright}p{2cm}>{\raggedright}p{2cm}>{\raggedright}p{2cm}>{\raggedright}p{2cm}>{\raggedright}p{2cm}}
\hline 
Considered period of time for the average & Average rank of support scenarios & Average rank of nonsupport scenarios & Average discharge of support scenarios ($10^{6}\,\text{m\textthreesuperior/day}$) & Average discharge of nonsupport scenarios ($10^{6}\,\text{m\textthreesuperior/day}$)\tabularnewline
 &  &  &  & \tabularnewline
\hline 
Year & $12.111$ & $11.077$ & $0.981$ & $0.857$\tabularnewline
Wet season & $12.889$ & $10.538$ & $1.338$ & $1.097$\tabularnewline
Dry season & $14.222$ & $9.615$ & $0.473$ & $0.517$\tabularnewline
Driest six months & $11.444$ & $11.539$ & $0.842$ & $0.938$\tabularnewline
Driest three months & $13.778$ & $9.923$ & $0.618$ & $0.651$\tabularnewline
Driest month & $14.333$ & $9.538$ & $0.418$ & $0.488$\tabularnewline
\hline 
\end{tabular}

\end{table}

\subsection{Computational time analysis}

As shown in Table \ref{tab:Computational-efficiency}, the computations
are very efficient, while they use the whole set of available river
discharges (twenty-three hydrological years). The most time-consuming
computations took a little more than half an hour. It corresponds
to the MPC algorithm being used with a weekly time step and mixing
scenario generation, each iteration having to deal with $23^{3}=12,\!167$
scenarios. Hence, both formalisms lead to efficient computation, without
specific implementation care, which is a positive aspect in terms
of scaleability of the method. 

The time ratio between the merging and mixing scenario-generation
techniques is about one thousand, while the size of the optimisation
problems (measured by the number of variables and constraints) grows
by a factor of approximately five hundred: this difference is explained
by the fact that the optimisation solver uses algorithms that do not
have a linear complexity in the number of variables \cite{Vanderbei2014}.
So far, no specific mitigation has been implemented to lower the computational
cost, such as taking advantage of the decoupled nature of the scenarios
(optimisation for each scenario can be performed independently from
the others). 

\begin{table}
\caption{\label{tab:Computational-efficiency} Computational efficiency comparison
between the different optimisation models. Historical data spans over
23 years. \protect \\
Preprocessing time (such as reading input files) is not included in
this benchmark. These tests are performed on a current high-end laptop
(Intel i7-6700HQ, 4 cores, 2.6 GHz, 16 GB of RAM) with the solver
Gurobi 6.5.0 and the modelling language Julia 0.4.7 (under Windows
10, 64 bits). }

\begin{tabular}{lll}
\hline 
Model & Time (seconds) & Number of scenarios\tabularnewline
 &  & \tabularnewline
\hline 
Stochastic (merge) & $1.8$ & $23-1=22$\tabularnewline
Stochastic (mix) & $2140.2$ & $23^{3}=12,\!167$\tabularnewline
Robust ($95\%$) & $0.1340$ & $1$\tabularnewline
Robust ($99\%$) & $0.1346$ & $1$\tabularnewline
\hline 
\end{tabular}
\end{table}

\section{Conclusion \label{sec:Conclusion} }

This article compares two paradigms to take uncertainty into account
within mathematical optimisation techniques applied to rule-curve
derivation: one is similar to many existing tools (stochastic programming),
with scenario generation to help deal with limited data; the other
one is more synthetic and directly uses confidence levels (robust
programming). Our results show that the current operating guidelines
can be improved at some points. Moreover, the dam operator can implement
it in such a way that rule curves are regularly updated. 

Besides those practical issues, the proposed methodology is easy to
implement efficiently enough, while being based on sound mathematical
principles. Also, the mathematical models used in this article are
exploited to their utmost potential: there is no point in pursuing
the research to get better solutions to these models, as the global
optimality has been reached. Nevertheless, the models could still
be improved to get more detailed results. For example, the discharge
through the hydropower penstock and the bottom outlets both depend
on the hydraulic head; they are currently considered as constants.
Another weak point is the computation of inflow confidence intervals,
which is crude, and could be greatly enhanced. 

Furthermore, this approach is directly applicable to any kind of water
demand, such as controlling low flows, as long as there is no uncertainty
in this demand; otherwise, another stage of uncertainty modelling
is needed, applying the same techniques as developed in this article. 

The rule curves have been evaluated with feasibility-related criteria,
which are the most relevant in this case for the operator. Nevertheless,
other evaluation processes could bring more information about the
behaviour of each potential policy with respect to the other dam purposes,
such as hydropower like in \cite{Arunkumar2012}. 

The approach can also be extended to handle flooding, by defining
a \emph{maximum} rule curve for normal operations: this lets some
free space to store the excess water due to flood events. This could
be computed, season per season, by analysing the maximum level so
that the usual constraints are not violated (as done in Section \ref{sec:Methodology}),
with the flood event as input. This requires some flood detection
algorithm, such as the one presented in \cite{Klopstra1999}. 
\begin{acknowledgements}
The authors gratefully acknowledge the \emph{Service Public de Wallonie}
(SPW) for providing data on the case study. 
\end{acknowledgements}

\section*{\textemdash \textemdash \textemdash \textemdash \textemdash \textemdash \textemdash{}}

\bibliographystyle{plain}
\bibliography{bs}

\begin{thebibliography}{10}

\bibitem{Abrahart1998}
Robert~J. Abrahart and Linda See.
\newblock {Neural Network vs. ARMA Modelling: constructing benchmark case
  studies of river flow prediction}.
\newblock In {\em Proceedings of the 3rd International Conference on
  GeoComputation}, Bristol, 1998.

\bibitem{Adam2014}
Nicolas Adam, S{\'{e}}bastien Erpicum, Pierre Archambeau, Michel Pirotton, and
  Benjamin Dewals.
\newblock {Stochastic Modelling of Reservoir Sedimentation in a Semi-Arid
  Watershed}.
\newblock {\em Water Resources Management}, 29(3):785--800, 2014.

\bibitem{Ahmad2014}
Asmadi Ahmad, Ahmed El-Shafie, Siti Fatin~Mohd Razali, and Zawawi~Samba
  Mohamad.
\newblock {Reservoir Optimization in Water Resources: a Review}.
\newblock {\em Water Resources Management}, 28(11):3391--3405, 2014.

\bibitem{Arunkumar2012}
R~Arunkumar and V~Jothiprakash.
\newblock {Optimal Reservoir Operation for Hydropower Generation using
  Non-linear Programming Model}.
\newblock {\em Journal of The Institution of Engineers (India): Series A},
  93(2):111--120, 2012.

\bibitem{Bashiri-Atrabi2015}
Hamid Bashiri-Atrabi, Kourosh Qaderi, David~E Rheinheimer, and Erfaneh Sharifi.
\newblock {Application of Harmony Search Algorithm to Reservoir Operation
  Optimization}.
\newblock {\em Water Resources Management}, 29(15):5729--5748, 2015.

\bibitem{Becker2014}
Berhard P.~J. Becker, Tobias Schruff, and Dirk Schwanenberg.
\newblock {Modellierung von reaktiver Steuerung und Model Predictive Control}.
\newblock In {\em 37. Dresdner Wasserbaukolloquium 2014}, 2014.

\bibitem{Ben-Tal2009}
Aharon Ben-Tal, Laurent {El Ghaoui}, and Arkadi Nemirovski.
\newblock {\em {Robust optimization}}.
\newblock Princeton University Press, 2009.

\bibitem{Ben-Tal2002}
Aharon Ben-Tal and Arkadi Nemirovski.
\newblock {Robust optimization: methodology and applications}.
\newblock {\em Mathematical Programming}, 92(3):453--480, 2002.

\bibitem{Bezanson2014}
Jeff Bezanson, Alan Edelman, Stefan Karpinski, and Viral~B. Shah.
\newblock {Julia: A Fresh Approach to Numerical Computing}.
\newblock nov 2014.

\bibitem{Bieri2013}
M.~Bieri and A.J. Schleiss.
\newblock {Analysis of flood-reduction capacity of hydropower schemes in an
  Alpine catchment area by semidistributed conceptual modelling}.
\newblock {\em Journal of Flood Risk Management}, 6(3):169--185, sep 2013.

\bibitem{Birge2011}
John~R Birge and Fran{\c{c}}ois Louveaux.
\newblock {\em {Introduction to Stochastic Programming}}.
\newblock Springer Verlag, second edition, 2011.

\bibitem{Camnasio2011}
Erica Camnasio and Gianfranco Becciu.
\newblock {Evaluation of the Feasibility of Irrigation Storage in a Flood
  Detention Pond in an Agricultural Catchment in Northern Italy}.
\newblock {\em Water Resources Management}, 25(5):1489--1508, 2011.

\bibitem{Celeste2009}
Alcigeimes~B. Celeste and Max Billib.
\newblock {Evaluation of stochastic reservoir operation optimization models}.
\newblock {\em Advances in Water Resources}, 32(9):1429--1443, 2009.

\bibitem{Chang2005}
Fi-John Chang, Li~Chen, and Li-Chiu Chang.
\newblock {Optimizing the reservoir operating rule curves by genetic
  algorithms}.
\newblock {\em Hydrological Processes}, 19(11):2277--2289, jul 2005.

\bibitem{Cortes1995}
Corinna Cortes and Vladimir Vapnik.
\newblock {Support-vector networks}.
\newblock {\em Machine Learning}, 20(3):273--297, 1995.

\bibitem{Dunning2015}
Iain Dunning, Joey Huchette, and Miles Lubin.
\newblock {JuMP: A Modeling Language for Mathematical Optimization}.
\newblock {\em arXiv:1508.01982 [math.OC]}, aug 2015.

\bibitem{Finch2008}
Jon Finch and Ann Calver.
\newblock {Methods for the quantification of evaporation from lakes}.
\newblock Technical report, 2008.

\bibitem{Klopstra1999}
D.~Klopstra and N.~Vrisou van Eck.
\newblock {Methodiek voor vaststelling van de vorm van de maatgevende
  afvoergolf van de Maas bij Borgharen}.
\newblock 1999.

\bibitem{Kwon2005}
Wook~Hyun Kwon and Soo~He Han.
\newblock {\em {Receding Horizon Control: Model Predictive Control for State
  Models}}.
\newblock Springer-Verlag London, 1 edition, 2005.

\bibitem{Labadie2004}
John~W. Labadie.
\newblock {Optimal Operation of Multireservoir Systems: State-of-the-Art
  Review}.
\newblock {\em Journal of Water Resources Planning and Management},
  130(2):93--111, 2004.

\bibitem{Maier2014}
Holger~Robert Maier, Zoran Kapelan, Joseph Kasprzyk, J~Kollat, Loren~Shawn
  Matott, M~C Cunha, Graeme~Clyde Dandy, Matt~S. Gibbs, Edward~C. Keedwell,
  Angela Marchi, Avi Ostfeld, Dragan Savic, Dimitri~P. Solomatine, Jasper~A.
  Vrugt, Aaron~C. Zecchin, B~S Minsker, Emily~J. Barbour, G~Kuczera, Fayzul
  Pasha, Andrea Castelletti, Matteo Giuliani, and P~M Reed.
\newblock {Evolutionary algorithms and other metaheuristics in water resources:
  Current status, research challenges and future directions}.
\newblock {\em Environmental Modelling {\&} Software}, 62:271--299, 2014.

\bibitem{Pan2015}
Limeng Pan, Mashor Housh, Pan Liu, Ximing Cai, and Xin Chen.
\newblock {Robust stochastic optimization for reservoir operation}.
\newblock {\em Water Resources Research}, 51(1):409--429, 2015.

\bibitem{Schwanenberg2015}
D~Schwanenberg, B~P~J Becker, and M~Xu.
\newblock {The open real-time control (RTC)-Tools software framework for
  modeling RTC in water resources sytems}.
\newblock {\em Journal of Hydroinformatics}, 17(1):130 LP -- 148, jan 2015.

\bibitem{Shapiro2014}
Alexander Shapiro and Darinka Dentcheva.
\newblock {\em {Lectures on stochastic programming: modeling and theory}},
  volume~16.
\newblock SIAM, second edition, 2014.

\bibitem{Taghian2014}
Mehrdad Taghian, Dan Rosbjerg, Ali Haghighi, and Henrik Madsen.
\newblock {Optimization of Conventional Rule Curves Coupled with Hedging Rules
  for Reservoir Operation}.
\newblock {\em Journal of Water Resources Planning and Management},
  140(5):693--698, may 2014.

\bibitem{Talsma2013}
Jan Talsma, Simone Patzke, Berhard P.~J. Becker, Neeltje Goorden, Dirk
  Schwanenberg, and Geert Prinsen.
\newblock {Application of Model Predictive Control on Water Extractions in
  Scarcity Situations in the Netherlands}.
\newblock {\em Revista de Ingenier{\'{i}}a Innova}, 6:1--10, 2013.

\bibitem{Vanderbei2014}
Robert~J Vanderbei.
\newblock {\em {Linear Programming - Foundations and Extensions}}.
\newblock Linear Programming, 4th edition, 2014.

\bibitem{WRCR:WRCR3909}
William W-G. Yeh.
\newblock {Reservoir Management and Operations Models: A State-of-the-Art
  Review}.
\newblock {\em Water Resources Research}, 21(12):1797--1818, 1985.

\bibitem{Zhang2015}
Yanke Zhang, Zhiqiang Jiang, Changming Ji, and Ping Sun.
\newblock {Contrastive analysis of three parallel modes in multi-dimensional
  dynamic programming and its application in cascade reservoirs operation}.
\newblock {\em Journal of Hydrology}, 529, Part:22--34, 2015.

\bibitem{Zhao2014}
T~Zhao, J~Zhao, J~Lund, and D~Yang.
\newblock {Optimal Hedging Rules for Reservoir Flood Operation from Forecast
  Uncertainties}.
\newblock {\em Journal of Water Resources Planning and Management},
  Preview(2011), 2014.

\end{thebibliography}

\appendix

\section{Technical characteristics }

\begin{table}[H]
\caption{\label{tab:Technical-characteristics}Characteristics of the catchment,
reservoir, and dam. }

\begin{tabular}{>{\raggedright}p{5cm}l}
\hline 
 & Value for the Eupen reservoir\tabularnewline
 & \tabularnewline
\hline 
Capacity & $25\cdot10^{6}\,\text{hm}^{3}$\tabularnewline
Dam height & $66\,\text{m}$\tabularnewline
Natural rivers & Vesdre, Getzbach\tabularnewline
Natural drainage area & $6920\,\text{ha}$\tabularnewline
Diverted river & Helle\tabularnewline
Diverted drainage area & $3675\,\text{ha}$\tabularnewline
Reservoir level for drinking water & $343\,\text{m}$\tabularnewline
Reservoir level for chosen flood safety margin & $361\,\text{m}$\tabularnewline
Spillway crest & $362\,\text{m}$\tabularnewline
\hline 
\end{tabular}
\end{table}

\begin{table}[H]
\caption{\label{tab:Technical-characteristics-optimisation}Characteristics
of the catchment, reservoir, and dam that are used for the optimisation. }

\begin{tabular}{>{\raggedright}p{5cm}ll}
\hline 
 & Symbol & Value for the Eupen reservoir\tabularnewline
 &  & \tabularnewline
\hline 
Minimum stored volume for drinking water  & $\mathrm{minStorage}$ & $2.25\,\text{hm}^{3}$\tabularnewline
Maximum stored volume for flood safety margin & $\mathrm{maxStorage}$ & $22\,\text{hm}^{3}$\tabularnewline
Required amount of drinking water & $\mathrm{drinkingWater}$ & $55\,000\text{ m\textthreesuperior/day}$\tabularnewline
Minimum environmental flow from the dam & $\mathrm{environmentalFlow}$ & $0.22\,\text{m\textthreesuperior/s}$\tabularnewline
Maximum discharge through the hydropower penstock & $\mathrm{penstockHydropower}$ & $4.5\,\text{m\textthreesuperior/s}$\tabularnewline
Maximum discharge through the bottom outlets & $\mathrm{bottomOutlet}$ & $70\,\text{m\textthreesuperior/s}$\tabularnewline
Maximum deviation from the river Helle & $\mathrm{maxDischarge}_{\mathrm{Helle}}$ & $20\,\text{m\textthreesuperior/s}$\tabularnewline
Minimum environmental flow for river Helle & $\mathrm{environmentalFlow}_{\mathrm{Helle}}$ & $0.05\,\text{m\textthreesuperior/s}$\tabularnewline
\hline 
\end{tabular}
\end{table}

\section{Impact of the confidence level on the robust solutions}

\begin{figure}[H]
\begin{centering}
\includegraphics[scale=0.4]{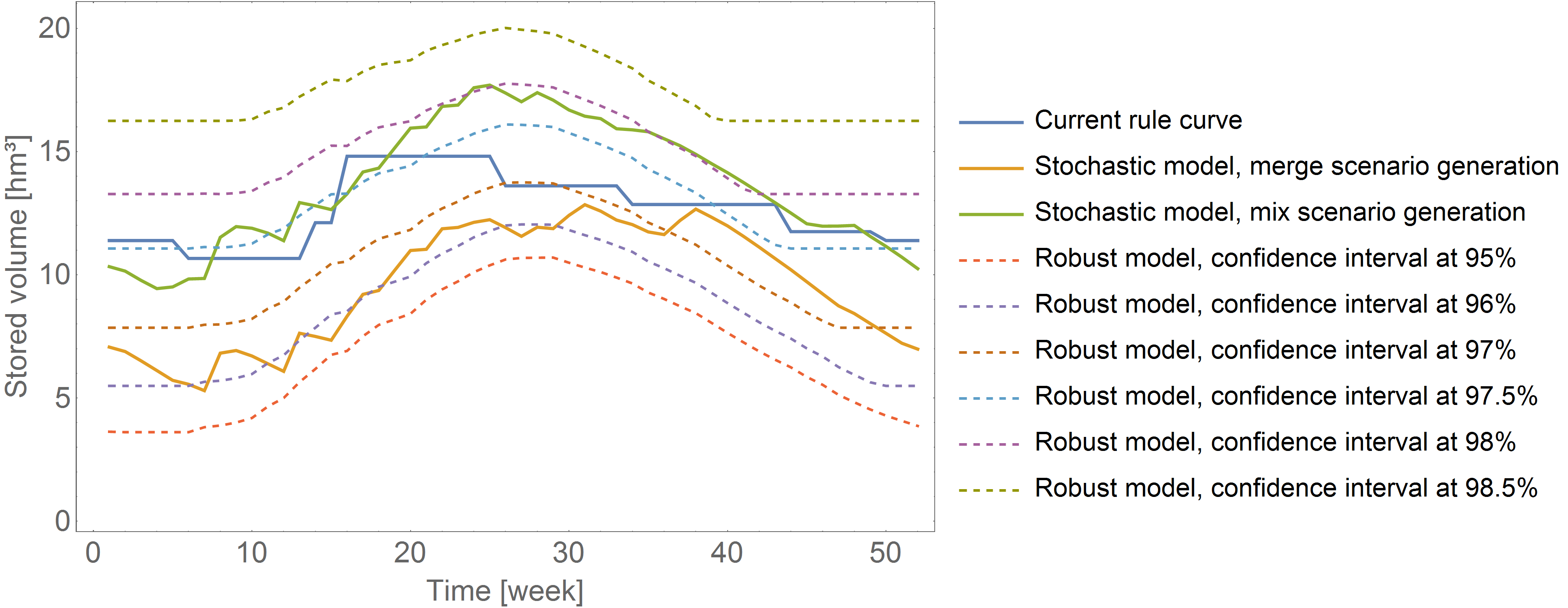}
\par\end{centering}
\caption{\label{fig:Compare-rule-curves-confidence-level}Impact of the confidence
level on the solutions.}
\end{figure}

\section{Determination of important years to perform the stochastic optimisation}

\begin{figure}[H]
\begin{centering}
\subfloat[\label{fig:Impact-scenario-on-solution-driest-three-months}The colours
correspond to the average discharge during \emph{the driest three-month}
\emph{period}. ]{\begin{centering}
\includegraphics[scale=0.3]{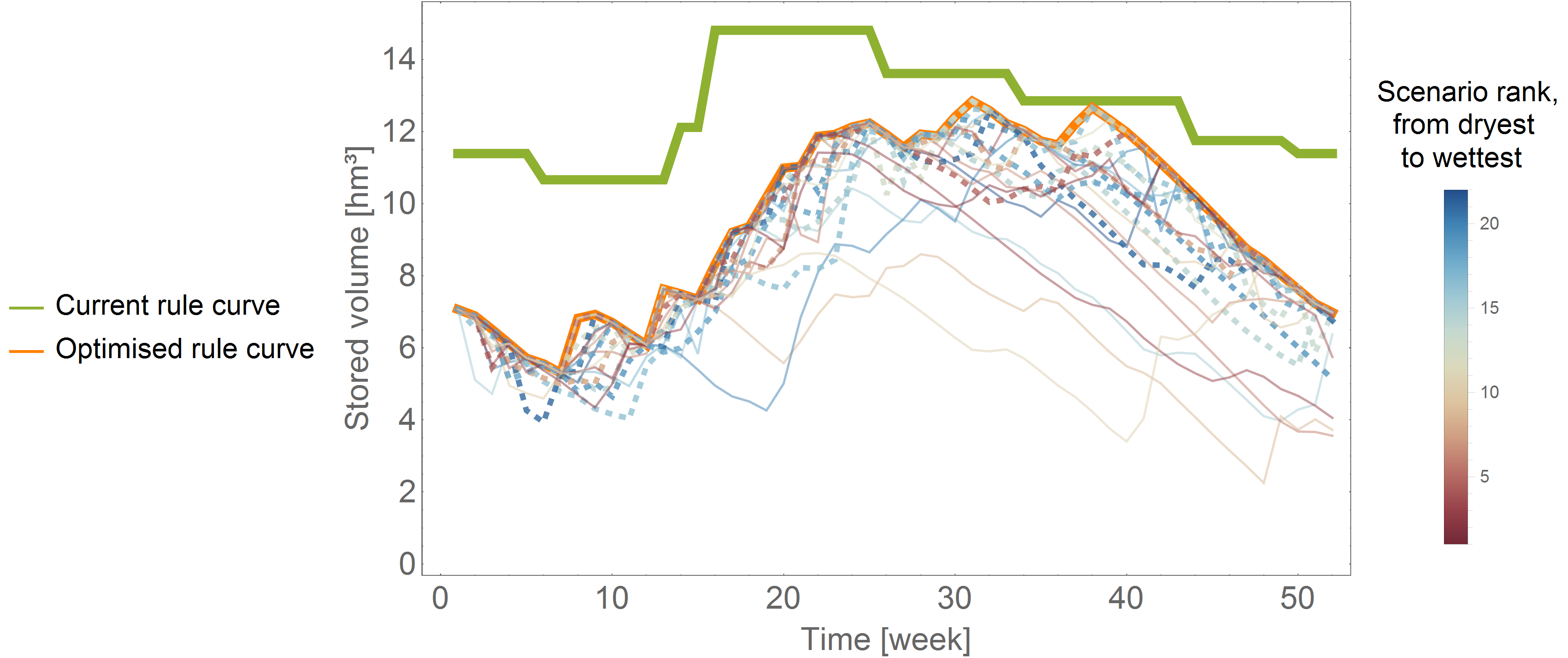}
\par\end{centering}
}
\par\end{centering}
\begin{centering}
\subfloat[\label{fig:Impact-scenario-on-solution-driest-six-months}The colours
correspond to the average discharge during \emph{the driest six-month}
\emph{period}. ]{\begin{centering}
\includegraphics[scale=0.3]{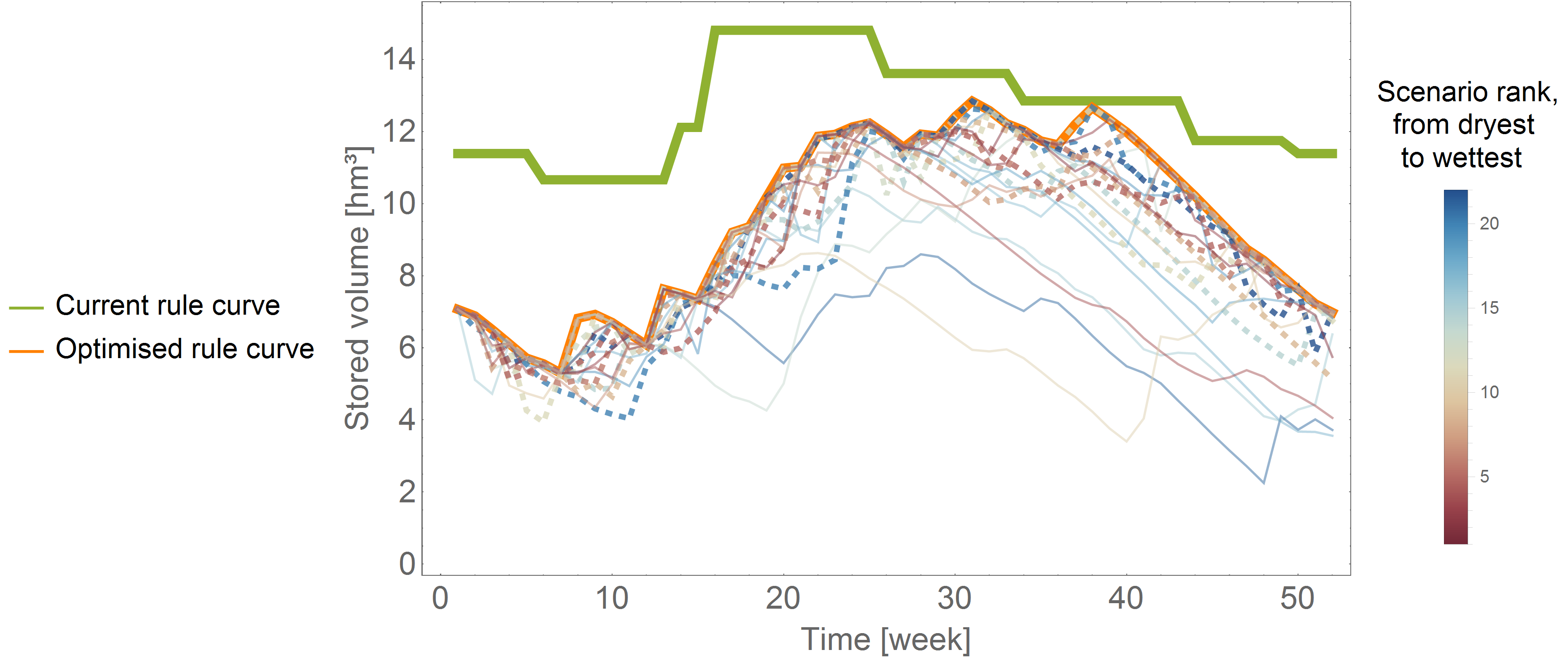}
\par\end{centering}
}
\par\end{centering}
\begin{centering}
\subfloat[\label{fig:Impact-scenario-on-solution-dry}The colours correspond
to the average discharge during \emph{the dry season} (May-September). ]{\begin{centering}
\includegraphics[scale=0.3]{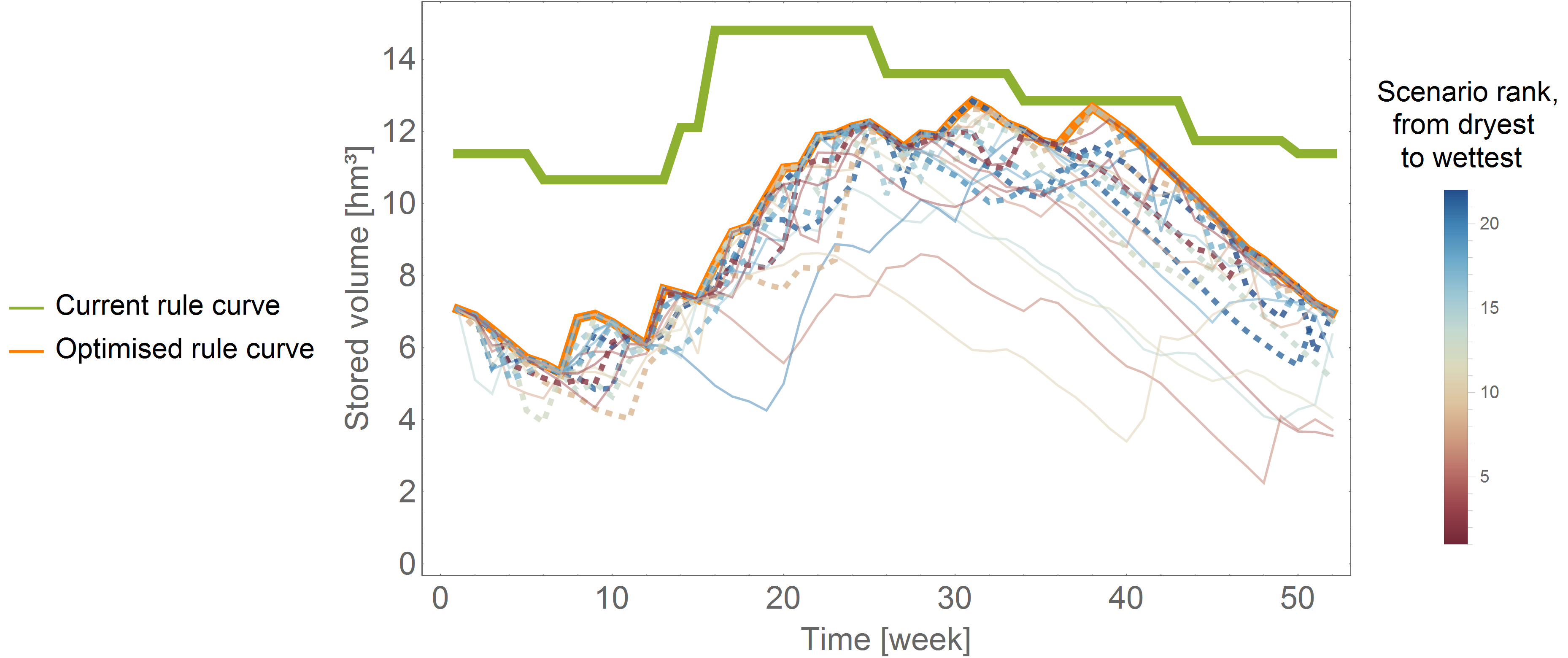}
\par\end{centering}
}
\par\end{centering}
\begin{centering}
\subfloat[\label{fig:Impact-scenario-on-solution-wet}The colours correspond
to the average discharge during \emph{the wet season} (October-April). ]{\begin{centering}
\includegraphics[scale=0.3]{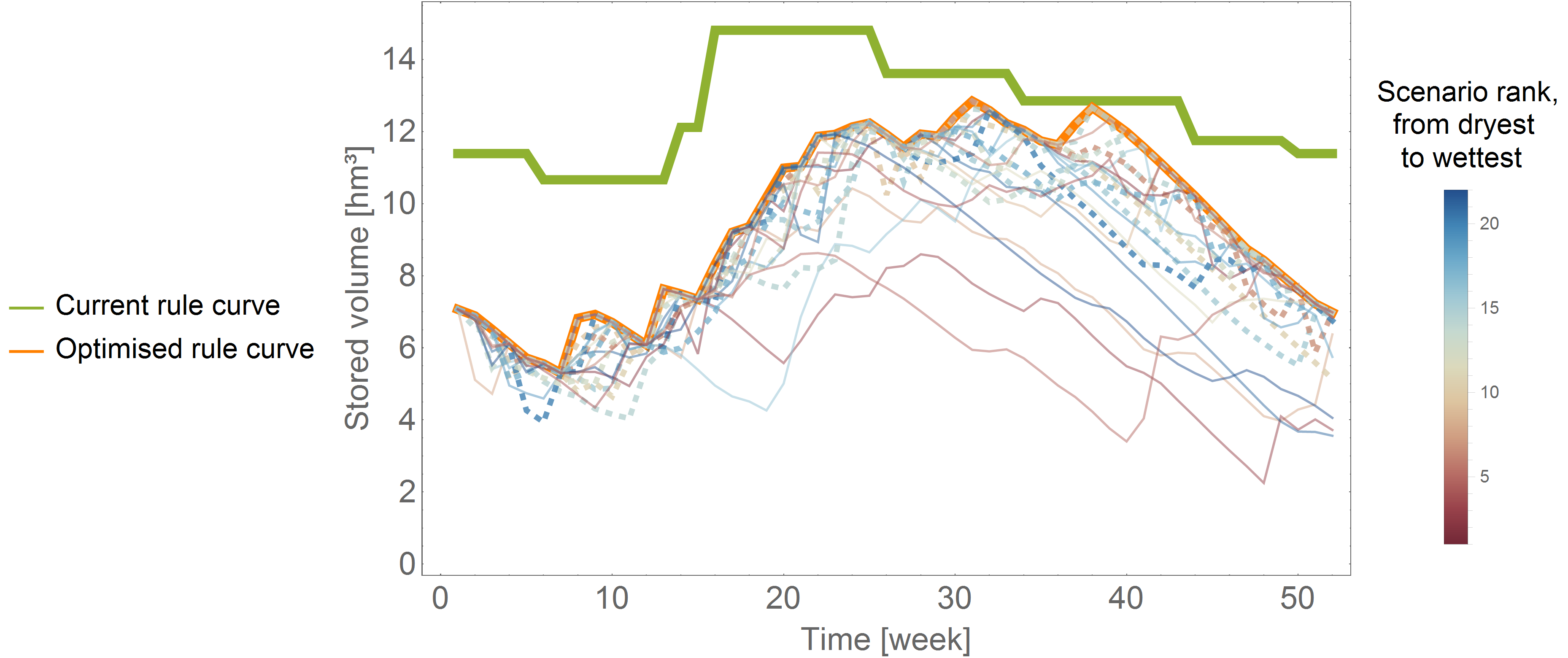}
\par\end{centering}
}
\par\end{centering}
\caption{\label{fig:Impact-scenario-on-solution-supplementary}Supplementary
analysis (refers to Figures \ref{fig:Impact-scenario-on-solution-yearly}
and \ref{fig:Impact-scenario-on-solution-driest-month}). }
\end{figure}

\end{document}